\documentclass[12pt,leqno]{article}
\usepackage{amsmath,amssymb,bbm}
\usepackage{amsthm}
\def\bm#1{\mathbbm{#1}}
\def\fn#1{\mathop{{\rm #1}\vphantom{\dim}}}

\makeatletter
\renewcommand{\section}{\@startsection{section}{1}{0mm}{12mm}{5mm}{\raggedright\bf\Large}}

\def\@citex[#1]#2{\if@filesw\immediate\write\@auxout{\string\citation{#2}}\fi
  \def\@citea{}\@cite{\@for\@citeb:=#2\do
    {\@citea\def\@citea{\@citesep}\@ifundefined
       {b@\@citeb}{{\bf ?}\@warning
       {Citation `\@citeb' on page \thepage \space undefined}}%
{\csname b@\@citeb\endcsname}}}{#1}}
\def\@citesep{; }
\makeatother

\newtheoremstyle{Kang}{}{}{\itshape}{}{\bf}{}{.5em}{}
\theoremstyle{Kang}
\newtheorem{theorem}{}[section]

\newtheoremstyle{Kremark}{}{}{}{}{\bf}{}{.5em}{}
\theoremstyle{Kremark}

\newtheorem{other}{}
\newenvironment{idef}[1]{\renewcommand{\theother}{#1}\begin{other}}{\end{other}}
\newenvironment{Case}[1]{\medskip {\it Case #1.}}{}

\title{Noether's problem for $p$-groups with a \\ cyclic subgroup of index $p^2$}
\author{Ming-chang Kang\\[2mm]
Department of Mathematics and \\Taida Institute of Mathematical Sciences\\
National Taiwan University\\
Taipei, Taiwan\\
E-mail: kang@math.ntu.edu.tw}
\date{}

\begin{document}

\maketitle

\footnote{\hspace*{-7.5mm}
Partially Supported by National Center for Theoretic Science
(Taipei office).}

\begin{abstract}
{\noindent\bf Abstract}
Let $K$ be any field and $G$ be a finite group.
Let $G$ act on the rational function field $K(x_g:g\in G)$ by $K$-automorphisms defined by $g\cdot x_h=x_{gh}$
for any $g,h\in G$.
Noether's problem asks whether the fixed field $K(G)=K(x_g:g\in G)^G$ is rational
(=purely transcendental) over $K$.
We will prove that if $G$ is a non-abelian $p$-group of order $p^n$ ($n\ge 3$) containing a cyclic subgroup
of index $p^2$ and $K$ is any field containing a primitive $p^{n-2}$-th root of unity,
then $K(G)$ is rational over $K$.
As a corollary, if $G$ is a non-abelian $p$-group of order $p^3$ and $K$ is a field containing
a primitive $p$-th root of unity, then $K(G)$ is rational.
\end{abstract}

\section{Introduction}

Let $K$ be any field and $G$ be a finite group. Let $G$ act on the
rational function field $K(x_g:g\in G)$ by $K$-automorphisms such
that $g\cdot x_h=x_{gh}$ for any $g,h\in G$. Denote by $K(G)$ the
fixed field $K(x_g:g\in G)^G$. Noether's problem asks whether
$K(G)$ is rational (=purely transcendental) over $K$. It is
related to the inverse Galois problem, to the existence of generic
$G$-Galois extensions over $K$, and to the existence of versal
$G$-torsors over $K$-rational field extensions \cite[33.1,
p.86]{Sw,Sa1,GMS}.  Noether's problem for abelian groups was
studied extensively by Swan, Voskresenskii, Endo, Miyata and
Lenstra, etc. The reader is referred to Swan's paper for a survey
of this problem \cite{Sw}.

On the other hand, just a handful of results about Noether's
problem are obtained when the groups are not abelian. In this
article we will restrict our attention to Noether's problem for
non-abelian $p$-groups.

First we recall several known results of along this direction.

\begin{theorem}[Chu and Kang \cite{CK}] \label{t1.1}
Let $G$ be a non-abelian $p$-group of order $\le p^4$ and exponent $p^e$.
Assume that $K$ is any field such that either
{\rm (i)} $\fn{char}K=p>0$, or {\rm (ii)} $\fn{char}K\ne p$ and $K$ contains a primitive $p^e$-th root of unity.
Then $K(G)$ is rational over $K$.
\end{theorem}

\begin{theorem}[Kang \cite{Ka1}] \label{t1.2}
Let $G$ be a non-abelian metacyclic $p$-group of exponent $p^e$.
Assume that $K$ is any field such that either {\rm (i)} $\fn{char}K=p>0$,
or {\rm (ii)} $\fn{char}K\ne p$ and $K$ contains a primitive $p^e$-th root of unity.
The $K(G)$ is rational over $K$.
\end{theorem}

\begin{theorem}[Saltman \cite{Sa2}] \label{t1.3}
Let $K$ be any field with $\fn{char}K\ne p$ $($in particular,
$K$ may be any algebraically closed field with $\fn{char}K\ne p)$.
There exists a non-abelian $p$-group $G$ of order $p^9$ such that $K(G)$ is not rational over $K$.
\end{theorem}

\begin{theorem}[Bogomolov \cite{Bo}] \label{t1.4}
There exists a non-abelian $p$-group $G$ of order $p^6$ such that $\bm{C}(G)$ is not rational over $\bm{C}$.
\end{theorem}

\begin{theorem}[Chu, Hu, Kang and Prokhorov \cite{CHKP}] \label{t1.5}
Let $G$ be a non-abelian group of order $32$ and exponent $2^e$.
Assume that $K$ is a field satisfying that either {\rm (i)} $\fn{char}K=2$,
or {\rm (ii)} $\fn{char}K\ne 2$ and $K$ contains a primitive $2^e$-th root of unity.
Then $K(G)$ is rational over $K$.
\end{theorem}

\begin{theorem}[Chu, Hu, Kang and Kunyavskii \cite{CHKK}] \label{t1.6}
Let $G$ be a non-abelian group of order $64$ and $K$ be a quadratically closed field
$($in particular, $\fn{char}K\ne 2)$.
Denote by $B_0(G,\mu)$ the unramified Brauer group of $G$ over $K$
$($where $\mu$ is the multiplicative group of all roots of unity in $K\backslash\{0\})$,
and by $G(i)$ the $i$-th group in the database of GAP for groups of order $64$.
\begin{enumerate}
\item[$(1)$] The following statements are equivalent,
\begin{enumerate}
\item[{\rm (a)}] $B_0(G,\mu)\ne 0$;
\item[{\rm (b)}] $Z(G)\simeq C_2^2$, $[G,G]\simeq C_4\times C_2$, $G/[G,G]\simeq C_2^3$,
$G$ has no abelian subgroup of index $2$,
and $G$ has no faithful $4$-dimensional representation over $\bm{C}$;
\item[{\rm (c)}] $G$ is isomorphic to one of the nine groups $G(i)$ where
$i=149$, $150$, $151$, $170$, $171$, $172$, $177$, $178$, $182$.
\end{enumerate}
\item[$(2)$] If $B_0(G,\mu)\ne 0$, then $K(G)$ is not stably rational over $K$.
\item[$(3)$] If $B_0(G,\mu)=0$, then $K(G)$ is rational over $K$ except possibly for groups $G$
which is isomorphic to $G(i)$ with $241\le i\le 245$.
\end{enumerate}
\end{theorem}

\begin{theorem}[Hu and Kang \cite{HuK}] \label{t1.7}
Let $n\ge 3$ and $G$ be a non-abelian $p$-group of order $p^n$ such that $G$ contains a cyclic subgroup of index $p$.
Assume that $K$ is any field satisfying that either {\rm (i)} $\fn{char}K=p>0$,
or {\rm (ii)} $\fn{char}K\ne p$ and $K$ contains a primitive $p^{n-2}$-th root of unity.
Then $K(G)$ is rational over $K$.
\end{theorem}

The main result of this article is the following theorem,
which is a generalization of the above \ref{t1.7}.

\begin{theorem} \label{t1.8}
Let $n\ge 3$ and $G$ be a non-abelian $p$-group of order $p^n$ such that $G$ contains a cyclic subgroup
of index $p^2$.
Assume that $K$ is any field satisfying that either
{\rm (i)} $\fn{char}K=p>0$, or
{\rm (ii)} $\fn{char}K\ne p$ and $K$ contains a primitive $p^{n-2}$-th root of unity.
Then $K(G)$ is rational over $K$.
\end{theorem}

Using \ref{t1.7}, the proof of \ref{t1.8} consists of three
ingredients : (a) rationality criteria mentioned before and some
other ones to be summarized in the next section, (b)
classification of $p$-groups with a cyclic subgroup of index
$p^2$, which is due to Ninomiya (see Section 3), and (c) a case by
case study of the rationality problems for the groups in (b).
Although there are so many groups to be checked and a case by case
study looks formidable, the rationality problems of most of these
groups look rather similar. It turns out that there are only three
typical cases, i.e. Case 1 and Case 5 of Section 4 and Case 5 of
Section 5.

By the way, we remark that, if $K$ doesn't contain enough roots of
unity (e.g.\ $K=\bm{Q}$) and $G$ is a non-abelian $p$-group, the
rationality of $K(G)$ is known only for a few cases at present.
See \cite{CHK,Ka2,Ka4} and the references therein.

By \ref{t1.8}, it is possible to simplify the proof of \ref{t1.1}
as follows. Using \ref{t1.8}, to show that $K(G)$ is rational when
$G$ is a non-abelian group with order $p^3$ or $p^4$, it suffices
to consider the rationality problem of $K(G)$ where $G$ is a
non-abelian $p$-group of order $p^4$ and exponent $p$ (such that
$K$ is a field containing a primitive $p$-th root of unity). There
are only two non-isomorphic groups of this type, i.e.\ (VI) and
(VII) in \cite[Theorem 3.2]{CK}. The rationality of $K(G)$ for
these two groups can be proved by the same method as in \cite{CK}.

We organize this article as follows. Section 2 contains more
rationality criteria which will be used subsequently. In Section
3, we recall the classification of non-abelian $p$-groups with a
cyclic subgroup of index $p^2$ by Ninomiya \cite{Ni}, which was
reproved by Berkovich and Janko \cite{BJ1}. The proof of
\ref{t1.8} is given in Section 4 and Section 5.

\begin{idef}{Standing Notations.}
Throughout this article, $K(x_1,\ldots,x_n)$ or $K(x,y)$ will be
rational function fields over $K$. $\zeta_n$ denotes a primitive
$n$-th root of unity. Whenever we write $\fn{char}K \nmid n$, it
is understood that either $\fn{char}K=0$ or $\fn{char}K>0$ with
$\gcd\{n,\fn{char}K\}=1$. When we write $\zeta_n\in K$, it is
assumed tacitly that $\fn{char}K \nmid n$. A field extension $L$
of $K$ is called rational over $K$ (or $K$-rational, for short) if
$L\simeq K(x_1,\ldots,x_n)$ over $K$ for some integer $n$. $L$ is
stably rational over $K$ if $L(y_1,\ldots,y_m)$ is rational over
$K$ for some $y_1,\ldots,y_m$ which are algebraically independent
over $L$. Recall that $K(G)$ denotes $K(x_g:g\in G)^G$ where
$h\cdot x_g=x_{hg}$ for $h,g\in G$.

A group $G$ is called metacyclic, if $G$ can be generated by two
elements $\sigma$ and $\tau$, and one of them generates a normal
subgroup of $G$. $C_n$ denotes the cyclic group of order $n$. The
exponent of a finite group $G$ is $\fn{lcm}\{\fn{ord}(g):g\in G\}$
where $\fn{ord}(g)$ is the order of $g$.

If $G$ is a finite group acting on a rational function field $K(x_1,\ldots,x_n)$ by $K$-automorphisms,
the actions of $G$ are called purely monomial actions if, for any $\sigma\in G$, any $1\le j\le n$,
$\sigma\cdot x_j=\prod_{1\le i\le n}x_i^{a_{ij}}$ where $a_{ij}\in\bm{Z}$;
similarly, the actions of $G$ are called monomial actions if, for any $\sigma\in G$,
any $1\le j\le n$, $\sigma\cdot x_j=\lambda_j(\sigma)\cdot\prod_{1\le i\le n}x_i^{a_{ij}}$ where
$a_{ij}\in\bm{Z}$ and $\lambda_j(\sigma)\in K\backslash\{0\}$.
All the groups in this article are finite groups.
\end{idef}

\section{Preliminaries}

In this section we recall several results which will be used in the proof of \ref{t1.8}.

\begin{theorem}[Kuniyoshi \cite{Ku}] \label{t2.1}
Let $K$ be a field with $\fn{char}K=p>0$ and $G$ be a $p$-group.
Then $K(G)$ is rational over $K$.
\end{theorem}

\begin{theorem}[{\cite[Theorem 1]{HK}}] \label{t2.2}
Let $G$ be a finite group acting on $L(x_1,\ldots,x_n)$,
the rational function field of $n$ variables over a field $L$.
Suppose that
\begin{enumerate}
\item[{\rm (i)}] for any $\sigma\in G$, $\sigma(L)\subset L$;
\item[{\rm (ii)}] the restriction of the action of $G$ to $L$ is faithful;
\item[{\rm (iii)}] for any $\sigma\in G$,
\[
\begin{pmatrix} \sigma(x_1) \\ \sigma(x_2) \\ \vdots \\ \sigma(x_n) \end{pmatrix}
=A(\sigma)\cdot \begin{pmatrix} x_1 \\ x_2 \\ \vdots \\ x_n \end{pmatrix}+B(\sigma)
\]
where $A(\sigma)\in GL_n(L)$ and $B(\sigma)$ is an $n\times 1$ matrix over $L$.
\end{enumerate}

Then there exist elements $z_1,\ldots,z_n\in L(x_1,\ldots,x_n)$
such that $L(x_1,\ldots,x_n)=L(z_1,\ldots,z_n)$ and
$\sigma(z_i)=z_i$ for any $\sigma\in G$, any $1\le i\le n$.
\end{theorem}

\begin{theorem}[{\cite[Theorem 3.1]{AHK}}] \label{t2.3}
Let $L$ be any field, $L(x)$ the rational function field of one
variable over $L$, and $G$ a finite group acting on $L(x)$.
Suppose that, for any $\sigma\in G$, $\sigma(L)\subset L$ and
$\sigma(x)=a_\sigma\cdot x+b_\sigma$ where $a_\sigma,b_\sigma\in
L$ and $a_\sigma\ne 0$. Then $L(x)^G=L^G(f)$ for some polynomial
$f\in L[x]$. In fact, if \ $m=\min \{\deg g(x) : g(x)\in L[x]^G
\setminus L \}$, any polynomial $f\in L[x]^G$ with $\deg f=m$
satisfies the property $L(x)^G=L^G(f)$.
\end{theorem}

\begin{theorem}[{\cite[Theorem 1.9]{KP}}] \label{t2.4}
Let $K$ be any field, $G_1$ and $G_2$ be two finite groups.
If both $K(G_1)$ and $K(G_2)$ are rational over $K$, so is $K(G_1\times G_2)$.
\end{theorem}

\begin{theorem}[\cite{Ha}] \label{t2.5}
Let $G$ be a finite group acting on the rational function field
$K(x,y)$ by monomial $K$-automorphisms. Then $K(x,y)^G$ is
rational over $K$.
\end{theorem}

\begin{theorem}[Fischer {\cite[Theorem 6.1; KP, Corollary 1.5]{Sw}}] \label{t2.6}
Let $G$ be a finite abelian group of exponent $e$,
and let $K$ be a field containing a primitive $e$-th root of unity.
For any linear representation $G\to GL(V)$ over $K$,
the fixed field $K(V)^G$ is rational over $K$.
\end{theorem}

\begin{theorem}[{\cite[Theorem 1.4]{Ka3}}] \label{t2.7}
Let $K$ be a field and $G$ be a finite group. Assume that {\rm
(i)} $G$ contains an abelian normal subgroup $H$ so that $G/H$ is
cyclic of order $n$, {\rm (ii)} $\bm{Z}[\zeta_n]$ is a unique
factorization domain, and {\rm (iii)} $\zeta_e\in K$ where $e$ is
the exponent of $G$. If $G\to GL(V)$ is any finite-dimensional
linear representation of $G$ over $K$, then $K(V)^G$ is rational
over $K$.
\end{theorem}

\section{Ninomiya's Theorem}

Let $n\ge 3$ and $p$ be a prime number. A complete list of
non-abelian $p$-groups of order $p^n$ containing a cyclic subgroup
of index $p$ was given by Burnside early in 1911 (see, for
examples, \cite[p.107; HuK, Theorem 1.9; Ni, p.1]{Su}). The
classification of non-abelian $p$-groups of order $p^n$ containing
a cyclic subgroup of index $p^2$ was completed rather late. This
problem was investigated by Burnside, G.\ A.\ Miller, etc.\ (see
\cite[Remark 3]{Ni}). The classification problem was solved by
Yasushi Ninomiya in 1994 \cite{Ni}. M. Kumar and L. Vermani,
apparently ignorant of Ninomiya's paper, provides a partial list
of these groups in \cite{KV}. Unfortunately their list contained
some mistakes, which were detected in \cite[p.31--32]{FN}. A
different proof of Ninomiya's Theorem was given by Berkovich and
Janko \cite[Section 11; BJ2, Section 74]{BJ1}. Now we state
Ninomiya's Theorem.

\begin{theorem}[Ninomiya {\cite[Theorem 1]{Ni}}] \label{t3.1}
Let $n\ge 3$ and $p$ be an odd prime number.
The finite non-abelian $p$-groups of order $p^n$ which have a cyclic subgroup of index $p^2$,
but haven't a cyclic subgroup of order $p$ are of the following types:
\begin{enumerate}
\item[{\rm (I)}] $n\ge 3$ \\
$G_1=\langle \sigma,\tau,
\lambda:\sigma^{p^{n-2}}=\tau^p=\lambda^p=1,\sigma\lambda=\lambda\sigma,
\tau\lambda=\lambda\tau,\tau^{-1}\sigma\tau=\sigma\lambda\rangle$.
\item[{\rm (II)}] $n\ge 4$ \\
$\begin{array}{@{}r@{\;}l}
G_2= & \langle \sigma,\tau:\sigma^{p^{n-2}}=\tau^{p^2}=1,\tau^{-1}\sigma\tau=\sigma^{1+p^{n-3}}\rangle, \\
G_3= & \langle \sigma,\tau,\lambda:\sigma^{p^{n-2}}=\tau^p=\lambda^p=1,\sigma\lambda=\lambda\sigma,
\tau\lambda=\lambda\tau,\tau^{-1}\sigma\tau=\sigma^{1+p^{n-3}}\rangle, \\
G_4= & \langle \sigma,\tau,\lambda:\sigma^{p^{n-2}}=\tau^p=\lambda^p=1,\sigma\tau=\tau\sigma,
\sigma\lambda=\lambda\sigma,\lambda^{-1}\tau\lambda=\sigma^{p^{n-3}}\tau\rangle, \\
G_5= & \langle \sigma,\tau,\lambda:\sigma^{p^{n-2}}=\tau^p=\lambda^p=1,\sigma\tau=\tau\sigma,
\lambda^{-1}\sigma\lambda=\sigma\tau,\lambda^{-1}\tau\lambda=\sigma^{p^{n-3}}\tau\rangle, \\
G_6= &\langle \sigma,\tau,\lambda: \sigma^{p^{n-2}}=\tau^p=\lambda^p=1,\sigma\tau=\tau\sigma,
\lambda^{-1}\sigma\lambda=\sigma\tau,\lambda^{-1}\tau\lambda=\sigma^{a\cdot p^{n-3}}\tau\rangle \\
& \mbox{where $\bar{a}\in \bm{Z}/p\bm{Z} \backslash \{\bar{0}\}$ is a quadratic non-residue}, \\
G_7= & \langle \sigma,\tau,\lambda: \sigma^{p^{n-2}}=\tau^p=\lambda^p=1,
\tau^{-1}\sigma\tau=\sigma^{1+p^{n-3}},\lambda^{-1}\sigma\lambda=\sigma\tau,\tau\lambda=\lambda\tau\rangle.
\end{array}$
\item[{\rm (III)}] $n\ge 5$ \\
$G_8=\langle \sigma,\tau:\sigma^{p^{n-2}}=\tau^{p^2}=1,\tau^{-1}\sigma\tau=\sigma^{1+p^{n-4}}\rangle$, \\
$G_9=\langle \sigma,\tau:\sigma^{p^{n-2}}=\tau^{p^2}=1,\sigma^{-1}\tau\sigma=\tau^{1+p}\rangle$.
\item[{\rm (IV)}] $n\ge 6$ \\
$G_{10}=\langle \sigma,\tau: \sigma^{p^{n-2}}=1, \sigma^{p^{n-3}}=\tau^{p^2},\sigma^{-1}\tau\sigma=\tau^{1-p}\rangle$.
\item[{\rm (V)}] $n=4$ and $p=3$ \\
$G_{11}=\langle \sigma,\tau,\lambda:\sigma^9=\tau^3=1,\sigma^3=\lambda^3,\sigma\tau=\tau\sigma,
\lambda^{-1}\sigma\lambda=\sigma\tau,\lambda^{-1}\tau\lambda=\sigma^6\tau\rangle$.
\end{enumerate}
\end{theorem}

\begin{theorem}[Ninomiya {\cite[Theorem 2]{Ni}}] \label{t3.2}
Let $n\ge 4$.
The finite non-abelian groups of order $2^n$ which have a cyclic subgroup of index $4$,
but haven't a cyclic subgroup of index $2$ are of the following types:
\begin{enumerate}
\item[{\rm (I)}] $n\ge 4$ \\
$G_1=\langle\sigma,\tau:\sigma^{2^{n-2}}=\tau^4=1,\tau^{-1}\sigma\tau=\sigma^{1+2^{n-3}}\rangle$, \\
$G_2=\langle\sigma,\tau,\lambda: \sigma^{2^{n-2}}=\lambda^2=1,\sigma^{2^{n-3}}=\tau^2,
\tau^{-1}\sigma\tau=\sigma^{-1},\sigma\lambda=\lambda\sigma,\tau\lambda=\lambda\tau\rangle$, \\
$G_3=\langle\sigma,\tau,\lambda: \sigma^{2^{n-2}}=\tau^2=\lambda^2=1,
\tau^{-1}\sigma\tau=\sigma^{-1},\sigma\lambda=\lambda\sigma,\tau\lambda=\lambda\tau\rangle$, \\
$G_4=\langle\sigma,\tau,\lambda:\sigma^{2^{n-2}}=\tau^2=\lambda^2=1,\sigma\tau=\tau\sigma,
\sigma\lambda=\lambda\sigma,\lambda^{-1}\tau\lambda=\sigma^{2^{n-3}}\tau\rangle$, \\
$G_5=\langle \sigma,\tau,\lambda: \sigma^{2^{n-2}}=\tau^2=\lambda^2=1,\sigma\tau=\tau\sigma,
\lambda^{-1}\sigma\lambda=\sigma\tau,\tau\lambda=\lambda\tau\rangle$.
\item[{\rm (II)}] $n\ge 5$ \\
$\begin{array}{@{}l@{\;}c@{\;}l}
G_6 &=& \langle\sigma,\tau:\sigma^{2^{n-2}}=\tau^4=1,\tau^{-1}\sigma\tau=\sigma^{-1}\rangle, \\
G_7 &=& \langle\sigma,\tau:\sigma^{2^{n-2}}=\tau^4=1,\tau^{-1}\sigma\tau=\sigma^{-1+2^{n-3}}\rangle, \\
G_8 &=& \langle\sigma,\tau:\sigma^{2^{n-2}}=1,\sigma^{2^{n-3}}=\tau^4,\tau^{-1}\sigma\tau=\sigma^{-1}\rangle, \\
G_9 &=& \langle\sigma,\tau:\sigma^{2^{n-2}}=\tau^4=1,\sigma^{-1}\tau\sigma=\tau^{-1}\rangle, \\
G_{10} &=& \langle\sigma,\tau,\lambda:\sigma^{2^{n-2}}=\tau^2=\lambda^2=1,\tau^{-1}\sigma\tau=\sigma^{1+2^{n-3}},
\sigma\lambda=\lambda\sigma,\tau\lambda=\lambda\tau\rangle, \\
G_{11} &=& \langle\sigma,\tau,\lambda:\sigma^{2^{n-2}}=\tau^2=\lambda^2=1,\tau^{-1}\sigma\tau=\sigma^{-1+2^{n-3}},
\sigma\lambda=\lambda\sigma,\tau\lambda=\lambda\tau\rangle, \\
G_{12} &=& \langle\sigma,\tau,\lambda:\sigma^{2^{n-2}}=\tau^2=\lambda^2=1,\sigma\tau=\tau\sigma,
\lambda^{-1}\sigma\lambda=\sigma^{-1},\lambda^{-1}\tau\lambda=\sigma^{2^{n-3}}\tau\rangle, \\
G_{13} &=& \langle\sigma,\tau,\lambda:\sigma^{2^{n-2}}=\tau^2=\lambda^2=1,\sigma\tau=\tau\sigma,
\lambda^{-1}\sigma\lambda=\sigma^{-1}\tau,\tau\lambda=\lambda\tau\rangle, \\
G_{14} &=& \langle\sigma,\tau,\lambda:\sigma^{2^{n-2}}=\tau^2=1,\sigma^{2^{n-3}}=\lambda^2,\sigma\tau=\tau\sigma,
\lambda^{-1}\sigma\lambda=\sigma^{-1}\tau,\tau\lambda=\lambda\tau\rangle, \\
G_{15} &=& \langle\sigma,\tau,\lambda:\sigma^{2^{n-2}}=\tau^2=\lambda^2=1,\tau^{-1}\sigma\tau=\sigma^{1+2^{n-3}},
\lambda^{-1}\sigma\lambda=\sigma^{-1+2^{n-3}},\tau\lambda=\lambda\tau\rangle, \\
G_{16} &=& \langle\sigma,\tau,\lambda:\sigma^{2^{n-2}}=\tau^2=\lambda^2=1,\tau^{-1}\sigma\tau=\sigma^{1+2^{n-3}},
\lambda^{-1}\sigma\lambda=\sigma^{-1+2^{n-3}}, \\
&& ~\lambda^{-1}\tau\lambda=\sigma^{2^{n-3}}\tau\rangle, \\
G_{17} &=& \langle\sigma,\tau,\lambda:\sigma^{2^{n-2}}=\tau^2=\lambda^2=1,\tau^{-1}\sigma\tau=\sigma^{1+2^{n-3}},
\lambda^{-1}\sigma\lambda=\sigma\tau,\tau\lambda=\lambda\tau\rangle, \\
G_{18} &=& \langle\sigma,\tau,\lambda:\sigma^{2^{n-2}}=\tau^2=1,\lambda^2=\tau,
\tau^{-1}\sigma\tau=\sigma^{1+2^{n-3}},\lambda^{-1}\sigma\lambda=\sigma^{-1}\tau\rangle.
\end{array}$
\item[{\rm (III)}] $n\ge 6$ \\
$\begin{array}{@{}r@{\;}l}
G_{19}=& \langle \sigma,\tau:\sigma^{2^{n-2}}=\tau^4=1,\tau^{-1}\sigma\tau=\sigma^{1+2^{n-4}}\rangle, \\
G_{20}=& \langle \sigma,\tau:\sigma^{2^{n-2}}=\tau^4=1,\tau^{-1}\sigma\tau=\sigma^{-1+2^{n-4}}\rangle, \\
G_{21}=& \langle \sigma,\tau:\sigma^{2^{n-2}}=1,\sigma^{2^{n-3}}=\tau^4,\tau^{-1}\sigma\tau=\tau^{-1}\rangle, \\
G_{22}=& \langle\sigma,\tau,\lambda:\sigma^{2^{n-2}}=\tau^2=\lambda^2=1,\sigma\tau=\tau\sigma,
\lambda^{-1}\sigma\lambda=\sigma^{1+2^{n-4}}\tau,\lambda^{-1}\tau\lambda=\sigma^{2^{n-3}}\tau\rangle, \\
G_{23}=& \langle\sigma,\tau,\lambda:\sigma^{2^{n-2}}=\tau^2=\lambda^2=1,\sigma\tau=\tau\sigma,
\lambda^{-1}\sigma\lambda=\sigma^{-1+2^{n-4}}\tau,\lambda^{-1}\tau\lambda=\sigma^{2^{n-3}}\tau\rangle, \\
G_{24}=& \langle\sigma,\tau,\lambda:\sigma^{2^{n-2}}=\tau^2=\lambda^2=1,\tau^{-1}\sigma\tau=\sigma^{1+2^{n-3}},
\lambda^{-1}\sigma\lambda=\sigma^{-1+2^{n-4}}\tau,\tau\lambda=\lambda\tau\rangle, \\
G_{25}=& \langle\sigma,\tau,\lambda:\sigma^{2^{n-2}}=\tau^2=1,\sigma^{2^{n-3}}=\lambda^2,
\tau^{-1}\sigma\tau=\sigma^{1+2^{n-3}},\lambda^{-1}\sigma\lambda=\sigma^{-1+2^{n-4}}\tau, \\
& ~\tau\lambda=\lambda\tau\rangle,
\end{array}$
\item[{\rm (IV)}] $n=5$ \\
$\begin{array}{@{}r@{\;}l} G_{26}=&
\langle\sigma,\tau,\lambda:\sigma^8=\tau^2=1,\sigma^4=\lambda^2,
\tau^{-1}\sigma\tau=\sigma^5,\lambda^{-1}\sigma\lambda=\sigma\tau,
\tau\lambda=\lambda\tau\rangle.
\end{array}$

\end{enumerate}
\end{theorem}

\section{Proof of \ref{t1.8} when \boldmath{$p\ge 3$}}

In this section we will prove \ref{t1.8} when $p$ is an odd prime number.

If $\fn{char}K=p>0$, apply \ref{t2.1}.
Thus $K(G)$ is rational over $K$.

From now on till the end of this section, we assume that
$\fn{char}K\ne p$ and $K$ contains a primitive $p^{n-2}$-th root
of unity where $G$ is a $p$-group of order $p^n$ with $n\ge 3$.

Throughout this section, we will denote by $\zeta=\zeta_{p^{n-2}}$
for a primitive $p^{n-2}$-th root of unity.

Suppose that $G$ contains a cyclic subgroup of index $p$. Then
$K(G)$ is rational over $K$ by \ref{t1.7}. Thus we may consider
only those groups $G$ which have no cyclic subgroup of index $p$,
i.e.\ $G$ is one of the 11 groups listed in \ref{t3.1}.

We explain the general strategy of our proof. Let $V$ be a
$K$-vector space whose dual space $V^*$ is defined as
$V^*=\bigoplus_{g\in G} K\cdot x(g)$ where $G$ acts on $V^*$ by
$h\cdot x(g)=x(hg)$ for any $h,g\in G$. Thus $K(V)^G=K(x(g):g\in
G)^G=K(G)$. We will find a faithful subspace $W=\bigoplus_{1\le
i\le k} K\cdot y_i$ of $V^*$. By \ref{t2.2}, $K(G)$ is rational
over $K(y_1,\ldots,y_k)^G$. In particular, if
$K(y_1,\ldots,y_k)^G$ is rational over $K$, so is $K(G)$ over $K$.
As we will see, this faithful subspace $W$ is constructed as an
induced representation of certain 2-dimensional (or 3-dimensional)
representation of some abelian subgroup of $G$. We will illustrate
this idea in Step 1 of Case 1 in the following proof of
\ref{t1.8}.

Now we begin to prove \ref{t1.8} for $p\ge 3$.

\begin{Case}{1} $G=G_1$ where $G_1$ is the group in \ref{t3.1}. \end{Case}

Step 1. Recall that $\zeta=\zeta_{p^{n-2}}$ and
$V^*=\bigoplus_{g\in G} K\cdot x(g)$ on which $G$ acts by the
regular representation.

Define $\omega=\zeta^{p^{n-3}}$.
Thus $\omega$ is a primitive $p$-th root of unity.

Define $X_1,X_2\in V^*$ be
\[
X_1=\sum_{0\le j\le p^{n-2}-1} x(\sigma^j),~~~ X_2=\sum_{0\le j\le p-1} x(\lambda^j).
\]

Note that $\sigma\cdot X_1=X_1$ and $\lambda\cdot X_2=X_2$.

Define $Y_1,Y_2\in V^*$ by
\[
Y_1=\sum_{0\le j\le p-1} \omega^{-j}\lambda^j\cdot X_1, ~~~
Y_2=\sum_{0\le j\le p^{n-2}-1} \zeta^{-j}\sigma^j\cdot X_2.
\]

It follows that
\begin{align*}
\sigma :~& Y_1\mapsto Y_1, ~ Y_2\mapsto \zeta Y_2,  \\
\lambda :~& Y_1\mapsto \omega Y_1, ~ Y_2 \mapsto Y_2.
\end{align*}

Thus $K\cdot Y_1+K \cdot Y_2$ is a representation space of the subgroup $\langle \sigma,\lambda\rangle$.

Define $x_i=\tau^i\cdot Y_1$, $y_i=\tau^i\cdot Y_2$ for $0\le i\le
p-1$. It is easy to verify that, for $0\le i\le p-1$,
\begin{align*}
\sigma :~ & x_i\mapsto \omega^i x_i, ~ y_i \mapsto \zeta y_i, \\
\tau :~ & x_0\mapsto x_1 \mapsto \cdots \mapsto x_{p-1} \mapsto x_0, \\
& y_0 \mapsto y_1 \mapsto \cdots \mapsto y_{p-1} \mapsto y_0, \\
\lambda :~ & x_i \mapsto \omega x_i, ~ y_i\mapsto y_i.
\end{align*}

We find that $Y=\bigl(\bigoplus_{0\le i\le p-1} K\cdot
x_i\bigr)\oplus\bigl(\bigoplus_{0\le i\le p-1}K\cdot y_i\bigr)$ is
a faithful $G$-subspace of $V^*$. Thus, by \ref{t2.2}, it suffices
to show that $K(x_i,y_i:0\le i\le p-1)^G$ is rational over $K$.

\bigskip
Step 2. For $1\le i\le p-1$, define $u_i=x_i/x_{i-1}$ and
$v_i=y_i/y_{i-1}$. Thus $K(x_i,y_i: 0\le i\le
p-1)=K(x_0,y_0,u_i,v_i:1\le i\le p-1)$ and, for every $\rho \in
G$,
\[
\rho\cdot x_0 \in K(u_i,v_i:1\le i\le p-1)\cdot x_0, ~~~ \rho \cdot y_0\in K(u_i,v_i:1\le i\le p-1)\cdot y_0,
\]
while the subfield $K(u_i,v_i:1\le i\le p-1)$ is invariant (as a
subfield) by the action of $G$, i.e.
\begin{align*}
\sigma:~ & u_i \mapsto \omega u_i, ~ v_i \mapsto v_i, \\
\lambda:~ & u_i\mapsto u_i, ~ v_i\mapsto v_i, \\
\tau:~ & u_1\mapsto u_2\mapsto \cdots \mapsto u_{p-1}\mapsto (u_1u_2\cdots u_{p-1})^{-1}, \\
& v_1 \mapsto v_2 \mapsto\cdots \mapsto v_{p-1}\mapsto (v_1v_2\cdots v_{p-1})^{-1}.
\end{align*}

Apply \ref{t2.3}.
We find that, if $K(u_i,v_i:1\le i\le p-1)^G$ is rational over $K$,
so is $K(x_i,y_i: 0\le i\le p-1)^G$ over $K$.
It remains to show that $K(u_i,v_i: 1\le i\le p-1)^G$ is rational over $K$.

Since $\lambda$ acts trivially on $K(u_i,v_i: 1 \le i\le p-1)$,
we find that $K(u_i,v_i:1\le i\le p-1)^G=K(u_i,v_i: 1\le i\le p-1)^{\langle \sigma,\tau\rangle}$.

\bigskip
Step 3. We will linearize the action of $\tau$ on
$v_1,\ldots,v_{p-1}$.

Define $t_0=1+v_1+v_1v_2+v_1v_2v_3+\cdots+v_1v_2\cdots v_{p-1}$,
$t_1=1/t_0$, $t_i=v_1v_2\cdots v_{i-1}/ t_0$ for $2\le i\le p$.
Note that $\sum_{1\le i\le p} t_i=1$, $K(v_i:1 \le i\le p-1)=K(t_i:1\le i\le p-1)$ and
\[
\tau: t_0\mapsto t_0/v_1, ~ t_1\mapsto t_2\mapsto\cdots\mapsto t_{p-1}\mapsto t_p
=1-t_1-t_2-\cdots- t_{p-1}.
\]

Thus $K(u_i,v_i:1\le i\le p-1)=K(u_i,t_i:1\le i\le p-1)$.

Define $T_i=t_i-(1/p)$ for $1\le i\le p-1$.
Then $\tau: T_1\mapsto T_2\mapsto\cdots \mapsto T_{p-1}\mapsto -T_1-\cdots- T_{p-1}$.

\bigskip
Step 4. Write $L=K(u_i:1\le i\le p-1)$ and consider $L(T_i: 1\le
i\le p-1)^{\langle \sigma,\tau\rangle}$. Note that the group
$\langle \sigma,\tau\rangle$ acts on the field $L$ as $\langle
\sigma,\tau\rangle/ \langle \sigma^p \rangle$ and is faithful on
$L$. Thus we may apply \ref{t2.2} to $L(T_i: 1\le i\le
p-1)^{\langle \sigma,\tau\rangle}$. It remains to show that
$K(u_i: 1\le i\le p-1)^{\langle \sigma,\tau\rangle}$ is rational
over $K$.

Define $z_1=u_1^p$, $z_i=u_i/u_{i-1}$ for $2\le i\le p-1$. Then
$K(u_i:1\le i\le p-1)^{\langle \sigma\rangle}=K(z_i:1\le i\le
p-1)$ and the action of $\tau$ is given by
\begin{align*}
\tau :~& z_1 \mapsto z_1z_2^p, \\
& z_2\mapsto z_3\mapsto\cdots \mapsto z_{p-1}\mapsto (z_1z_2^{p-1}z_3^{p-2}\cdots z_{p-1}^2)^{-1}
\mapsto z_1z_2^{p-2}z_3^{p-3}\cdots z_{p-2}^2z_{p-1}\mapsto z_2.
\end{align*}

Define $s_1=z_2$, $s_i=\tau^{i-1}\cdot z_2$ for $2\le i\le p-1$.
Then $K(z_i: 1\le i\le p-1)=K(s_i: 1\le i\le p-1)$ and
\[
\tau: s_1\mapsto s_2\mapsto\cdots\mapsto s_{p-1}\mapsto (s_1s_2\cdots s_{p-1})^{-1} \mapsto s_1.
\]

The action of $\tau$ can be linearized as in Step 3.
Thus $K(s_i: 1\le i\le p-1)^{\langle \tau\rangle}$ is rational over $K$ by \ref{t2.6}.
Done.

\bigskip
\begin{Case}{2} $G=G_2$. \end{Case}

$G$ is a metacyclic group.
Apply \ref{t1.2}.
We find that $K(G)$ is rational over $K$.

\bigskip
\begin{Case}{3} $G=G_3$. \end{Case}

Define $H=\langle \sigma, \tau\rangle$. Then $G\simeq H\times
C_p$. $K(H)$ is rational over $K$ by \ref{t1.2} (alternatively, by
\ref{t1.7}). $K(C_p)$ is rational over $K$ by \ref{t2.6}. Thus
$K(G)$ is rational over $K$ by \ref{t2.4}.

\bigskip
\begin{Case}{4} $G=G_4$. \end{Case}

By the same method as in Step 1 of Case 1, for the abelian
subgroup $\langle \sigma,\tau\rangle$, choose $Y_1,Y_2\in
V^*=\bigoplus_{g\in G} K\cdot x(g)$ such that $\sigma\cdot
Y_1=Y_1$, $\sigma\cdot Y_2=\zeta Y_2$, $\tau\cdot Y_1=\omega Y_1$,
$\tau\cdot Y_2=Y_2$ where $\zeta=\zeta_{p^{n-2}}$ and
$\omega=\zeta^{p^{n-3}}$.

Define $x_i=\lambda^i\cdot Y_1$, $y_i=\lambda^i \cdot Y_2$ for $0\le i\le p-1$.
It follows that, for $0\le i\le p-1$, we have
\begin{align*}
\sigma:~& x_i\mapsto x_i, ~ y_i\mapsto \zeta y_i, \\
\tau:~& x_i\mapsto \omega x_i, ~ y_i\mapsto \omega^i y_i, \\
\lambda:~& x_0\mapsto x_1\mapsto\cdots \mapsto x_{p-1}\mapsto x_0, \\
& y_0\mapsto y_1 \mapsto\cdots\mapsto y_{p-1}\mapsto y_0.
\end{align*}

It suffices to show that $K(x_i,y_i:0\le i\le p-1)^G$ is rational over $K$.
The proof is almost the same as in Case 1.
Define $u_i=x_i/x_{i-1}$, $v_i=y_i/y_{i-1}$ for $1\le i\le p-1$.
We have
\begin{align*}
\sigma:~& u_i\mapsto u_i, ~ v_i\mapsto v_i, \\
\tau:~& u_i\mapsto u_i, ~ v_i\mapsto \omega v_i, \\
\lambda:~& u_1\mapsto u_2\mapsto\cdots\mapsto u_{p-1}\mapsto (u_1\cdots u_{p-1})^{-1}, \\
& v_1 \mapsto v_2 \mapsto\cdots\mapsto v_{p-1} \mapsto (v_1\cdots v_{p-1})^{-1}.
\end{align*}

Compare with the situation in Case 1. It is not difficult to show
that $K(u_i,v_i: 1\le i\le p-1)^G$ is rational over $K$.

\bigskip
\begin{Case}{5} $G=G_5$. \end{Case}

Step 1.
For the abelian subgroup $\langle \sigma,\tau\rangle$,
find $Y_1$ and $Y_2$ by the same way as in Case 4.

Define $x_i$, $y_i$ where $0\le i\le p-1$ by the same formulae as
in Case 4. Note that, for $0\le i\le p-1$, we have
\begin{align*}
\sigma:~& x_i \mapsto \omega^i x_i, ~ y_i\mapsto \zeta\omega^{i\choose 2} y_i, \\
\tau:~& x_i \mapsto \omega x_i, ~ y_i\mapsto \omega^i y_i, \\
\lambda:~& x_0\mapsto x_1\mapsto \cdots \mapsto x_{p-1}\mapsto x_0, \\
& y_0 \mapsto y_1 \mapsto \cdots \mapsto y_{p-1} \mapsto y_0,
\end{align*}
where $\zeta=\zeta_{p^{n-2}}$ and $\omega=\zeta^{p^{n-3}}$.

Define $u_i=x_i/x_{i-1}$, $v_i=y_i/y_{i-1}$ for $1\le i\le p-1$.
It suffices to show that $K(u_i,v_i:1\le i\le p-1)^G$ is rational over $K$.
Note that
\begin{equation}
\begin{aligned}
\sigma:~& u_i\mapsto \omega u_i, ~ v_i\mapsto \omega^{i-1}v_i, \\
\tau:~& u_i\mapsto u_i, ~ v_i\mapsto\omega v_i, \\
\lambda:~& u_1\mapsto u_2 \mapsto\cdots\mapsto u_{p-1}\mapsto (u_1\cdots u_{p-1})^{-1}, \\
& v_1\mapsto v_2 \mapsto \cdots \mapsto v_{p-1}\mapsto (v_1\cdots v_{p-1})^{-1}.
\end{aligned} \label{eq1}
\end{equation}

It follows that $K(u_i,v_i: 1\le i\le p-1)^{\langle \tau \rangle}=K(u_i,V_i:1\le i\le p-1)$ where
$V_1=v_1^p$ and $V_i=v_i/v_{i-1}$ for $2\le i\le p-1$.

Note that $\sigma:V_1\mapsto V_1$, $V_i\mapsto \omega V_i$ for $2\le i\le p-1$.
Moreover, $K(u_i,V_i: 1\le i\le p-1)^{\langle\sigma\rangle}=K(z_i,w_i:1\le i\le p-1)$ where
$z_1=u_1^p$, $w_1=V_1$, $z_i=u_i/u_{i-1}$, $w_i=V_i/u_i$ for $2\le i\le p-1$.

The action of $\lambda$ is given by
\begin{equation}
\begin{aligned}
\lambda:~& z_1 \mapsto z_1z_2^p, ~
z_2\mapsto z_3\mapsto\cdots\mapsto z_{p-1}\mapsto (z_1z_2^{p-1}z_3^{p-2}\cdots z_{p-1}^2)^{-1}, \\
& w_1\mapsto z_1z_2^pw_1w_2^p, ~
w_2\mapsto w_3\mapsto\cdots\mapsto w_{p-1}\mapsto A\cdot(w_1w_2^{p-1}w_3^{p-2}\cdots w_{p-1}^2)^{-1},
\end{aligned} \label{eq2}
\end{equation}
where $A$ is some monomial in $z_1, z_2,\ldots,z_{p-1}$.

We will ``linearize" the above action.

\bigskip
Step 2.
We write the additive version of the multiplication action of $\lambda$ in Formula \eqref{eq2},
i.e.\ consider the $\bm{Z}[\pi]$-module $M=\bigl(\bigoplus_{1\le i\le p-1} \bm{Z}\cdot z_i\bigr) \oplus %
\bigl(\bigoplus_{1\le i\le p-1}\bm{Z}\cdot w_i\bigr)$\vspace*{3pt} corresponding to \eqref{eq2} where
$\pi=\langle \lambda\rangle$.
Thus $\lambda$ acts on the $\bm{Z}$-base $z_i$, $w_i$ ($1\le i\le p-1$) as follows,
\begin{align*}
\lambda:~& z_1\mapsto z_1+pz_2, \\
& z_2\mapsto z_3\mapsto\cdots\mapsto z_{p-1}\mapsto -z_1-(p-1)z_2-(p-2)z_3-\cdots -2z_{p-1}, \\
& w_1\mapsto w_1+pw_2+z_1+pz_2, \\
& w_2\mapsto w_3\mapsto\cdots\mapsto w_{p-1}\mapsto -w_1-(p-1)w_2-(p-2)w_3-\cdots -2w_{p-1}+B
\end{align*}
where $B\in \bigoplus_{1\le i\le p-1}\bm{Z}\cdot z_i$ \vspace*{3pt}
(in fact, $B=\log A$ when interpreted suitably).

Define $M_1=\bigoplus_{1\le i\le p-1}\bm{Z}\cdot z_i$, which is a
$\bm{Z}[\pi]$-submodule of $M$. Define $M_2=M/M_1$. \vspace*{3pt}

It follows that we have a short exact sequence of $\bm{Z}[\pi]$-modules
\begin{equation}
0 \to M_1 \to M \to M_2 \to 0. \label{eq3}
\end{equation}

It is easy to see that $M_1\simeq M_2$ as $\bm{Z}[\pi]$-modules.

By Step 4 of Case 1, $M_1$ is isomorphic to the $\bm{Z}[\pi]$-module $N=\bigoplus_{1\le i\le p-1}\bm{Z}\cdot s_i$
where $s_1=z_2$, $s_i=\lambda^{i-1}\cdot z_2$ for $2\le i\le p-1$, and
\[
\lambda: s_1\mapsto s_2\mapsto\cdots\mapsto s_{p-1}\mapsto -s_1-s_2-\cdots -s_{p-1}\mapsto s_1.
\]

Let $\Phi_p(T)\in \bm{Z}[T]$ be the $p$-th cyclotomic polynomial.
Since $\bm{Z}[\pi]\simeq \bm{Z}[T]/T^p-1$,
we find that $\bm{Z}[\pi]/\Phi_p(\lambda) \simeq \bm{Z}[T]/\Phi_p(T) \simeq \bm{Z}[\omega]$,
the ring of $p$-th cyclotomic integer.
Note that the $\bm{Z}[\pi]$-module $N$ can be regarded as a $\bm{Z}[\omega]$-module through the morphism
$\bm{Z}[\pi]\to \bm{Z}[\pi]/\Phi_p(\lambda)$.
When $N$ is regarded as a $\bm{Z}[\omega]$-module,
$N\simeq \bm{Z}[\omega]$ the rank-one free $\bm{Z}[\omega]$-module.

We claim that $M$ itself may be regarded as a $\bm{Z}[\omega]$-module,
i.e.\ $\Phi_p(\lambda)\cdot M=0$.

Return to the multiplicative notations in Step 1. Note that $z_i$
and $w_i$ (where $1\le i\le p-1$) are monomials in $u_i$ and $v_i$
(where $1\le i\le p-1$). The action of $\lambda$ on $u_i$, $v_i$
given in Formula \eqref{eq1} satisfies the relation \ $\prod_{0\le
j\le p-1} \lambda^j(u_i)=\prod_{0\le j\le p-1}
\lambda^j(v_i)=1$\vspace*{3pt} for any $1\le i\le p-1$. Using the
additive notations, we get $\Phi_p(\lambda)\cdot
u_i=\Phi_p(\lambda)\cdot v_i=0$ for $1\le i\le p-1$.
Hence $\Phi_p(\lambda)\cdot m=0$ for any $m\in M \subset \bigl(\bigoplus_{1\le i\le p-1} \bm{Z}\cdot u_i\bigr) %
\oplus \bigl(\bigoplus_{1\le i\le p-1}\bm{Z}\cdot v_i\bigr)$. \vspace*{3pt}

In particular, the short exact sequence of $\bm{Z}[\pi]$-modules in Formula \eqref{eq3} is a short exact sequence
of $\bm{Z}[\omega]$-modules.

Since $M_1\simeq M_2\simeq N$ is a free $\bm{Z}[\omega]$-module,
the short exact sequence in Formula \eqref{eq3} splits, i.e.\
$M\simeq M_1\oplus M_2$ as $\bm{Z}[\omega]$-modules, and so as
$\bm{Z}[\pi]$-modules also.

\medskip
We interpret the additive version of $M\simeq M_1\oplus M_2\simeq N^2$ in terms of
the multiplicative version as follows:
There exist $Z_i$, $W_i$ (where $1\le i\le p-1$) such that $Z_i$ (resp.\ $W_i$) are monomials in $z_j$
and $w_j$ for $1\le j\le p-1$ and $K(z_i,w_i: 1\le i\le p-1)=K(Z_i,W_i: 1\le i\le p-1)$;
moreover, $\lambda$ acts as
\begin{align*}
\lambda:~& Z_1 \mapsto Z_2 \mapsto\cdots\mapsto Z_{p-1}\mapsto (Z_1\cdots Z_{p-1})^{-1}, \\
& W_1\mapsto W_2\mapsto\cdots\mapsto W_{p-1}\mapsto (W_1\cdots W_{p-1})^{-1}.
\end{align*}

The above action can be linearized (see Step 3 of Case 1).
Thus $K(Z_i,W_i: 1\le i\le p-1)^{\langle \lambda \rangle}$ is rational over $K$ by \ref{t2.6}.
This finishes the proof.

\bigskip
\begin{Case}{6} $G=G_6$. \end{Case}

As in Case 5, for the abelian subgroup $\langle
\sigma,\tau\rangle$, find $Y_1$ and $Y_2$; and define $x_i,y_i\in
V^*=\bigoplus_{g\in G} K\cdot x(g)$ such that, for $1\le i\le
p-1$,
\begin{align*}
\sigma:~& x_i\mapsto \omega^i x_i, ~ y_i\mapsto \zeta\omega^{{i\choose 2}a} y_i, \\
\tau:~& x_i\mapsto \omega x_i, ~ y_i\mapsto \omega^{ia} y_i, \\
\lambda:~& x_0 \mapsto x_1 \mapsto\cdots \mapsto x_{p-1}\mapsto x_0, \\
& y_0 \mapsto y_1 \mapsto \cdots \mapsto y_{p-1} \mapsto y_0.
\end{align*}

We will prove that $K(x_i,y_i:0\le i\le p-1)^G$ is rational over $K$.
The proof is almost the same as in the previous Case 5.
For $1\le i\le p-1$, define $u_i=x_i/x_{i-1}$, $v_i=y_i/y_{i-1}$.
Then we get
\begin{align*}
\sigma:~& u_i\mapsto \omega u_i, ~ v_i\mapsto \omega^{(i-1)a} v_i, \\
\tau:~& u_i\mapsto u_i, ~ v_i\mapsto \omega^a v_i, \\
\lambda:~& u_1\mapsto u_2\mapsto\cdots\mapsto u_{p-1}\mapsto (u_1u_2\cdots u_{p-1})^{-1}, \\
& v_1\mapsto v_2\mapsto\cdots\mapsto v_{p-1}\mapsto (v_1v_2\cdots v_{p-1})^{-1}.
\end{align*}

Then $K(u_i,v_i: 1\le i\le p-1)^{\langle\tau\rangle}=K(u_i,V_i: 1\le i\le p-1)$ where $V_1=v_1^p$,
$V_i=v_i/v_{i-1}$ for $2\le i\le p-1$.
The action of $\sigma$ is given by
\[
\sigma: V_1\mapsto V_1, ~ V_i\mapsto \omega^a V_i
\]
for $2\le i\le p-1$.

Define $z_1=u_1^p$, $w_1=V_1$, $z_i=u_i/u_{i-1}$, $w_i=V_i/u_i^a$ for $2\le i\le p-1$.
We get $K(u_i,V_i: 1\le i\le p-1)^{\langle \sigma\rangle}=K(z_i,w_i:1\le i\le p-1)$.
The remaining proof is the same as in Case 5.

\bigskip
\begin{Case}{7} $G=G_7$. \end{Case}

As before, let $\zeta=\zeta_{p^{n-2}}$, $\xi=\zeta^p$, $\omega=\zeta^{p^{n-3}}$ and
find $Y_1,Y_2,Y_3\in\bigoplus_{g\in G} K\cdot x(g)$ such that
\begin{align*}
\sigma^p:~& Y_1 \mapsto \xi Y_1, ~ Y_2\mapsto Y_2, ~ Y_3\mapsto Y_3, \\
\tau:~& Y_1\mapsto Y_1, ~ Y_2\mapsto \omega Y_2, ~ Y_3 \mapsto Y_3, \\
\lambda:~& Y_1 \mapsto Y_1, ~ Y_2 \mapsto Y_2, ~ Y_3 \mapsto \omega Y_3.
\end{align*}

For $0\le i\le p-1$, define $x_i=\sigma^i Y_1$, $y_i=\sigma^i Y_2$, $z_i=\sigma^i Y_3$.
Note that
\begin{align*}
\sigma:~& x_0\mapsto x_1\mapsto\cdots\mapsto x_{p-1}\mapsto\xi x_0, \\
& y_0\mapsto y_1\mapsto\cdots\mapsto y_{p-1}\mapsto y_0, \\
& z_0\mapsto z_1\mapsto\cdots\mapsto z_{p-1}\mapsto z_0, \\
\sigma^p:~& x_i\mapsto \xi x_i, ~ y_i\mapsto y_i, ~ z_i\mapsto z_i, \\
\tau:~& x_i\mapsto \omega^{-i}x_i, ~ y_i\mapsto \omega y_i, ~ z_i \mapsto z_i, \\
\lambda:~& x_i\mapsto \omega^{i\choose 2}x_i, ~ y_i\mapsto \omega^{-i}y_i,~ z_i\mapsto \omega z_i.
\end{align*}

Define $u_i=x_i/x_{i-1}$, $v_i=y_i/y_{i-1}$, $w_i=z_i/z_{i-1}$ where $1\le i\le p-1$.
It remains to show that $K(u_i,v_i,w_i: 1\le i\le p-1)^G$ is rational over $K$.

Define $W_0=1+w_1+w_1w_2+\cdots+w_1w_2\cdots w_{p-1}$, $W_1=1/W_0$,
$W_i=w_1\cdots w_{i-1}/ W_0$ for $2\le i\le p-1$;
define $U_i=u_i/\zeta$ for $1\le i\le p-1$.
It is easy to check that $K(u_i,v_i,w_i: 1\le i\le p-1)=K(U_i,v_i,W_i:1\le i\le p-1)$ and
\begin{equation}
\begin{aligned}
\sigma:~& U_1\mapsto U_2\mapsto\cdots\mapsto U_{p-1}\mapsto (U_1\cdots U_{p-1})^{-1} \mapsto U_1, \\
& v_1\mapsto v_2\mapsto\cdots \mapsto v_{p-1} \mapsto (v_1\cdots v_{p-1})^{-1} \mapsto v_1, \\
& W_1\mapsto W_2\mapsto\cdots\mapsto W_{p-1}\mapsto 1-W_1-W_2-\cdots-W_{p-1}, \\
\tau:~& U_i \mapsto \omega^{-1}U_i, ~ v_i\mapsto v_i, ~ W_i\mapsto W_i, \\
\lambda:~& U_i \mapsto \omega^{i-1} U_i, ~ v_i \mapsto
\omega^{-1}v_i, ~ W_i\mapsto W_i.
\end{aligned} \label{eq4}
\end{equation}

By \ref{t2.2}, if $K(U_i,v_i: 1\le i\le p-1)^G$ is rational over $K$,
so is $K(U_i,v_i,W_i: 1\le i\le p-1)^G$ over $K$.
Thus it remains to show that $K(U_i,v_i:1\le i\le p-1)^G$ is rational over $K$.

Compare the actions of $\sigma$, $\tau$, $\lambda$ in Formula
\eqref{eq4} with those in Formula \eqref{eq1}. They look almost
the same. Use the same method in Case 5. We find that
$K(U_i,v_i:1\le i\le p-1)^{\langle \sigma,\tau,\lambda\rangle}$ is
rational over $K$.

\bigskip
\begin{Case}{8} $G=G_8, G_9,G_{10}$. \end{Case}

These groups are metacyclic $p$-groups.
Apply \ref{t1.2} to conclude that $K(G)$ is rational over $K$.

\bigskip
\begin{Case}{9} $G=G_{11}$. \end{Case}

This group is of order 81 and with exponent 9.
Apply \ref{t1.1}.
We find that $K(G)$ is rational over $K$. Done.

\section{Proof of \ref{t1.8} when \boldmath{$p=2$}}

The idea of the proof for this situation is the same as that in
Section 4.

Thanks to \ref{t2.1},
we may assume that $\fn{char}K\ne 2$ and $K$ contains $\zeta=\zeta_{2^{n-2}}$,
a primitive $2^{n-2}$-th root of unity.

If $G$ is a non-abelian group of order 8,
it is isomorphic to the dihedral group or the quaternion group.
Thus $K(G)$ is rational over $K$ by \cite[Proposition 2.6 and Theorem 2.7]{CHK}.

From now on, we assume $G$ is a non-abelian group of order $2^n$
with $n\ge 4$. Since \ref{t1.7} takes care of the case when $G$
has an element of order $2^{n-1}$, we may consider only the case
when $G$ has an element of order $2^{n-2}$, but hasn't elements of
order $2^{n-1}$. Hence we may use the classification of $G$
provided by \ref{t3.2}. Namely, we will consider only those 25
groups in \ref{t3.2}.

\begin{Case}{1} $G=G_1,G_6,G_7,G_8,G_9,G_{19},G_{20},G_{21}$ in \ref{t3.2}. \end{Case}

These groups are metacyclic groups. Apply \ref{t1.2}. Done.

\begin{Case}{2} $G=G_2,G_3,G_{10},G_{11},G_{12}$. \end{Case}

Each of these groups $G$ contains a subgroup $H$ such that $G\simeq H\times C_2$.
Moreover, $H$ has a cyclic subgroup of index 2.
For example, when $G=G_2$, take $H=\langle \sigma,\tau\rangle$.
Apply \ref{t1.7} and \ref{t2.4}.

\begin{Case}{3} $G=G_4,G_5,G_{13},G_{14},G_{22},G_{23}$. \end{Case}

Each of these groups $G$ contains an abelian normal subgroup $H$ of index 2.
Apply \ref{t2.7}.

\begin{Case}{4} $G=G_{26}$. \end{Case}

This group is of order 32 and with exponent 8. Apply \ref{t1.5}.

\begin{Case}{5} $G=G_{15}$. \end{Case}

Denote $\zeta=\zeta_{2^{n-2}}$. Define $\xi=\zeta^2$.

As in the proof of the previous section,
for the abelian subgroup $\langle \sigma^2,\tau\rangle$,
find $Y_1,Y_2 \in \bigoplus_{g\in G} K\cdot x(g)$ such that
\begin{align*}
\sigma^2:~& Y_1 \mapsto \xi Y_1, ~ Y_2\mapsto Y_2, \\
\tau:~& Y_1\mapsto Y_1, ~ Y_2 \mapsto -Y_2.
\end{align*}

Define $x_0=Y_1$, $x_1=\sigma Y_1$, $x_2=\lambda Y_1$,
$x_3=\lambda\sigma Y_1$, $y_0=Y_2$, $y_1=\sigma Y_2$, $y_2=\lambda
Y_2$, $y_3=\lambda\sigma Y_2$. It is easy to verify that
\begin{align*}
\sigma:~& x_0\mapsto x_1\mapsto \xi x_0, ~ x_2\mapsto -\xi^{-1}x_3, ~ x_3\mapsto -x_2, \\
& y_0\mapsto y_1\mapsto y_0, ~ y_2\mapsto y_3\mapsto y_2,\\
\tau:~& x_0\mapsto x_0, ~ x_1\mapsto -x_1,~ x_2\mapsto x_2, ~ x_3\mapsto -x_3, \\
& y_i\mapsto -y_i \mbox{ ~ for }0\le i\le 3, \\
\lambda:~& x_0\mapsto x_2\mapsto x_0, ~ x_1\mapsto x_3\mapsto x_1, ~ y_0\mapsto y_2\mapsto y_0,
~ y_1\mapsto y_3\mapsto y_1.
\end{align*}

It suffices to show that $K(x_i,y_i:0\le i\le 3)^G$ is rational
over $K$.

Since $G$ acts faithfully on $K(x_i:0\le i\le 3)$, we may apply
\ref{t2.2} to $K(x_i,y_i:0\le i\le 3)^G$. It follows that
$K(x_i,y_i:0\le i\le 3)^G$ is rational over $K(x_i:0\le i\le
3)^G$. It remains to show that $K(x_i:0\le i\le 3)^G$ is rational
over $K$.

Define $u_1=x_0/x_1$, $u_2=x_2/x_3$, $u_3=x_1/x_2$. Apply
\ref{t2.3} to $K(x_i:0\le i\le 3)=K(u_1,u_2,u_3,x_3)$. We find
that $K(u_1,u_2,u_3,x_3)^G=K(u_1,u_2,u_3)^G(w)$ for some element
$w$ fixed by $G$. It suffices to show that
$K(u_1,u_2,u_3)^G=K(u_1,u_2,u_3)^{\langle \sigma,\tau,\lambda
\rangle}$ is rational over $K$.

The action of $G$ is given as follows,
\begin{align*}
\sigma:~& u_1\mapsto \xi^{-1}/u_1,~ u_2\mapsto \xi^{-1}/u_2, ~ u_3\mapsto -\xi^2u_1u_2u_3, \\
\tau:~& u_i\mapsto -u_i \mbox{ ~ for } 1\le i\le 3, \\
\lambda:~& u_1 \mapsto u_2 \mapsto u_1, u_3 \mapsto 1/(u_1u_2u_3).
\end{align*}

In particular, $\sigma^2(u_1)=u_1$, $\sigma^2(u_2)=u_2$,
$\sigma^2(u_3)=\xi^2u_3$.

Define $u_4=u_3^{2^{n-4}}$. Then $K(u_1,u_2,u_3)^{\langle \sigma^2
\rangle}=K(u_1,u_2,u_4)$. Note that
$\sigma(u_4)=(u_1u_2)^{2^{n-4}}u_4$ (because $n \ge 5$),
$\tau(u_4)=u_4$, $\lambda(u_4)=1/((u_1u_2)^{2^{n-4}}u_4)$.

Define $z_3=(u_1u_2)^{2^{n-5}}u_4$. We find that
$\sigma(z_3)=-z_3$, $\tau(z_3)=z_3$, $\lambda(z_3)=1/z_3$.

Define $z_1=u_1u_2,z_2=u_1/u_2$. It follows that
$K(u_1,u_2,u_4)^{\langle \tau \rangle}=K(z_1,z_2,z_3)$. Moreover,
$\sigma(z_1)=\xi^{-2}/z_1,\sigma(z_2)=1/z_2$,
$\lambda(z_1)=z_1,\lambda(z_2)=1/z_2$.

Define $v=(1-z_2)/(1+z_2)$. Then $\sigma(v)=-v,\lambda(v)=-v$.
Apply \ref{t2.2} to $K(z_1,z_2,z_3)^{\langle \sigma,\lambda
\rangle}=K(z_1,z_3,v)^{\langle \sigma,\lambda \rangle}$. We find
that $K(z_1,z_3,v)^{\langle \sigma,\lambda \rangle}$ is rational
over $K(z_1,z_3)^{\langle \sigma,\lambda \rangle}$. Note that
$K(z_1,z_3)^{\langle \sigma,\lambda \rangle}$ is rational over $K$
by \ref{t2.5}. Hence the result.

\bigskip
\begin{Case}{6} $G=G_{16}$. \end{Case}

The proof is almost the same as the previous Case 5.
For the abelian subgroup $\langle \sigma^2,\tau\rangle$, find $Y_1$ and $Y_2$.
Define $x_i$, $y_i$ (where $0\le i\le 3$) by the same way and try to show that $K(x_i,y_i:0\le i\le 3)^G$ is rational.
The action of $G$ is given by
\begin{align*}
\sigma:~& x_0 \mapsto x_1\mapsto \xi x_0, ~ x_2\mapsto -\xi^{-1} x_3, ~ x_3\mapsto -x_2, \\
& y_0 \mapsto y_1 \mapsto y_0, ~ y_2 \mapsto y_3 \mapsto y_2, \\
\tau:~& x_0 \mapsto x_0, ~ x_1 \mapsto -x_1, ~ x_2\mapsto -x_2, ~ x_3 \mapsto x_3, \\
& y_i \mapsto -y_i \mbox{ ~ for } 0\le i\le 3, \\
\lambda:~& x_0\leftrightarrow x_2, ~ x_1\leftrightarrow x_3, ~ y_0\leftrightarrow y_2, ~ y_1\leftrightarrow y_3.
\end{align*}

It suffices to consider the rationality of $K(x_i:0\le i\le 3)^G$.
Define $u_1, u_2, u_3$ by the same formulae as in the previous
Case 5. It follows that the actions of $\sigma, \tau, \lambda$ on
$u_1, u_2, u_3$ are completely the same as in Case 5, except that
$\tau(u_3)=u_3$ in the present situation (in Case 5, we have
$\tau(u_3)=-u_3$). The proof is the same as Case 5 and the details
are omitted.

\bigskip
\begin{Case}{7} $G=G_{17}$. \end{Case}

We give two proofs for this case.

For the first proof, we may use the same method in Case 5 of this
section. For the abelian subgroup $\langle \sigma^2,\tau\rangle$,
find $Y_1$ and $Y_2$. Define $x_i$, $y_i$ for $0\le i\le 3$. Then
define $u_1, u_2, u_3$ by the same way as in Case 5. Note that
$\sigma^2(u_3)=-u_3$ in this case. Thus we define $u_4=u_3^2$ in
the present case (instead of defining $u_4=u_3^{2^{n-4}}$ as in
Case 5). Then define $z_3=u_1u_2u_4, z_1=u_1u_2,z_2=u_1/u_2$. We
find that $K(u_1, u_2, u_3)^{\langle \sigma^2, \tau
\rangle}=K(z_1, z_2, z_3)$. Moreover, $\sigma(z_1)=-\xi^{-2}/z_1$,
$\lambda(z_1)=z_1$. Since $-\xi^{-2}=\alpha^2$ where
$\alpha=\sqrt{-1}\xi^{-1} \in K$, we define $v=(\alpha
-z_1)/(\alpha +z_1)$. We find that $\sigma (v)=-v, \lambda (v)=v$.
Thus we may apply \ref{t2.5}. Done.

Alternatively, this group is the special case of $G_7$ in
\ref{t3.1} when $p=2$. Note that, in the Case 7 of Section 4, we
don't use anything whether $p$ is odd or even. Thus the proof is
still valid for this situation.

\bigskip
\begin{Case}{8} $G=G_{18}, G_{24}, G_{25}$. \end{Case}

Again the proof is almost the same,
but some modification should be carried out.
We illustrate the situation $G=G_{18}$ as follows.

Consider the abelian subgroup $\langle \sigma^2,\tau\rangle$ and
define $Y_1$, $Y_2$, $x_i$, $y_i$ (where $0\le i\le 3$) and $u_i$
where $1\le i\le 3$. By the same method as in Case 5, we can show
that $K(u_1, u_2, u_3)^{\langle \sigma^2, \tau \rangle}=K(z_1,
z_2, z_3)$.

Now consider the action of $\sigma, \lambda$ on $z_1,z_2,z_3$.
This time we will linearize the actions on $z_1$ by the same
formula as in the first proof of Case 7. The remaining proof is
the same as before.

The situation when $G=G_{24}$ or $G_{25}$ is the same as the
situation $G=G_{18}$. Done.

\newpage
\renewcommand{\refname}{\centering{References}}

\end{document}